\documentclass[12pt,reqno]{amsart}
\usepackage{amsmath, amssymb, amsfonts, xcolor}
\usepackage{hyperref}
\numberwithin{equation}{section}
\usepackage{amsthm}
\newtheorem{theorem}{Theorem}
\newtheorem{lemma}{Lemma}

\newtheorem{corollary}{Corollary}

\theoremstyle{definition}
\newtheorem{definition}{Definition}

\newtheorem{remark}{Remark}

\usepackage[a4paper, top = 1.2in, bottom = 1.2in, left = 1.2in, right = 1.2in]{geometry}
\usepackage{enumitem}
\newlist{steps}{enumerate}{1}
\setlist[steps, 1]{label = Step \arabic*:}

\newtheoremstyle{case}{}{}{}{}{}{:}{ }{}
\theoremstyle{case}

\usepackage{enumitem}
\numberwithin{equation}{section}

\newcommand{\N}{\mathbb{N}}
\newcommand{\Q}{\mathbb{Q}}

\newcommand{\Z}{\mathbb{Z}}

\newcommand{\e}{\mathbf{e}}
\newcommand{\f}{\mathbf{f}}

\newcommand{\mcO}{\mathcal{O}}

\newcommand{\mfP}{\mathfrak{P}}

\def\1{1\!\!1}

\title[Generalized Fruit Diophantine equation over $K$]{Generalized Fruit Diophantine equation over number fields}

\author[S. Sahoo]{Satyabrat Sahoo}
\address[S. Sahoo]{Yau Mathematical Sciences Center, Tsinghua University, Beijing 100084, China.}
\email{sahoos@mail.tsinghua.edu.cn}

\author[S. Laishram]{Shanta Laishram}
\address[S. Laishram]{Stat-Math Unit, Indian Statistical Institute 7, S. J. S. Sansanwal Marg, New Delhi, 110016, India.}
\email{shanta@isid.ac.in}

\date{\today}

\begin{document}
	\maketitle
	\begin{abstract}
			Let $K$ be a number field and $\mcO_K$ be the ring of integers of $K$. In this article, we study the solutions of the generalized fruit Diophantine equation $ax^d-y^2-z^2 +xyz-c=0$ over $K$, where $d \geq 3$ is an integer and $a,c\in \mcO_K\setminus \{0\}$. 
		Subsequently, we provide explicit values of square-free integers $t$ such that the equation $ax^d-y^2-z^2 +xyz-c=0$ has no solution $(x_0, y_0, z_0) \in \mcO_{\Q(\sqrt{t})}^3$ with $2 | x_0$, and demonstrate that the set of all such square-free integers $t$ with $t \geq 2$ has density exactly $\frac{1}{6}$.
		As an application, we construct infinitely many elliptic curves $E$ defined over number fields $K$ having no integral point $(x_0,y_0) \in \mcO_K^2$ with $2|x_0$.
	\end{abstract}
	
	\maketitle
	
	%    \tableofcontents
	\section{Introduction}
Diophantine equations are among the most active and fascinating areas in number theory. In \cite{W95}, Wiles proved the famous Fermat's Last Theorem using modularity and established that the Diophantine equation $x^n+y^n=z^n$ has no non-zero integer solutions for integers $n \geq 3$. Subsequently, significant progress has been achieved in the examination of the generalized Fermat equation $Ax^p+By^q=Cz^r$ over number fields (see \cite{D16}, \cite{FS15}, \cite{IKO20}, \cite{JS25},  \cite{KS23 Diophantine1}, \cite{KS23 Diophantine2}, \cite{M23}, \cite{S24 GFLT}, \cite{S24 GFDE}, \cite{KS24 GAFLT} for more details).

In \cite{LT08}, Luca and Togb\'e first studied the integer solutions of the Diophantine equation $x^3 + by+1-xyz = 0$ with a fixed integer $b$. Later in \cite{T09}, Togb\'e studied the integer solutions of the Diophantine equation $x^3 + by+4-xyz = 0$. In \cite{MS23}, Majumdar and Sury proved that the Diophantine equation $x^3-y^2-z^2 +xyz-5=0$ has no integer solutions and named this equation the \textit{fruit Diophantine equation}. As an application, in  \cite{MS23} Majumdar and Sury constructed infinitely many elliptic curves with no integral points. In \cite{VS22}, Vaishya and Sharma extended the work of \cite{MS23} to the Diophantine equation $ax^3-y^2-z^2 +xyz-b=0$ for integers $a,b$ with $a \equiv 1 \pmod {12}$ and $b=8a-3$, and constructed infinitely many elliptic curves with torsion-free Mordell--Weil group over $\Q$.
In \cite{PC23}, Prakash and Chakraborty generalized the work of \cite{VS22} to the generalized fruit Diophantine equation $ax^d-y^2-z^2 +xyz-b=0$ for integers $a,b$ with $a \equiv 1 \pmod {12}$ and $b=2^da-3$, where $d$ is an odd integer divisible by $3$, and constructed infinitely many hyperelliptic curves with torsion-free Mordell--Weil group over $\Q$.

In this article, we study the solutions of the generalized fruit Diophantine equation $ax^d-y^2-z^2 +xyz-c=0$ over number fields $K$, where $d \geq 3$ is an integer and $a,c\in \mcO_K\setminus \{0\}$. 

In Theorem~\ref{main result1 for GFDE}, we show that for any $a,b\in \mcO_K \setminus \{0\}$ and $c=2^db-3^r$ with integers $r\geq 2$ and $d \geq 3$ odd, the generalized fruit Diophantine equation $ax^d-y^2-z^2 +xyz-c=0$ has no solution $(x_0, y_0, z_0) \in \mcO_K^3$ with $2 | x_0$. As an application of this result, for almost all algebraic integers $\alpha \in \mcO_K$, we construct elliptic curves $E_\alpha$ defined over $K$ such that $E_\alpha$ has no integral point $(x_0,y_0) \in \mcO_K^2$ with $2|x_0$ (see Theorem~\ref{appl thm}). This generalizes the work of \cite{VS22} where they constructed elliptic curves $E_m$ defined over $\Q$ for integers $m$ such that $E_m$ has no point $(x_0,y_0) \in \Z^2$.

In Corollary~\ref{cor for quadratic field}, we provide explicit values of square-free integers $t$ such that the hypothesis of Theorem~\ref{main result1 for GFDE} holds over $K=\Q(\sqrt{t})$. Finally, in Theorem~\ref{thm for density}, we show that the set of square-free integers $t \geq 2$ such that the hypothesis of Theorem~\ref{main result1 for GFDE} holds over $K=\Q(\sqrt{t})$ has density exactly $\frac{1}{6}.$
%		\end{enumerate}

		\subsection{Structure of the article:}
The structure of this article is as follows. In Section~\ref{section for main result}, we state the main results of the article, namely, Theorems~\ref{main result1 for GFDE} and~\ref{thm for density}. In Section~\ref{steps to prove main results}, we prove Theorems~\ref{main result1 for GFDE} and~\ref{thm for density}. In Section~\ref{section for application}, we state and prove Theorem~\ref{appl thm}.		

\section{Main Results for the Generalized Fruit Diophantine Equation over $K$}
\label{section for main result}
Throughout this article, $K$ denotes a number field and $\mcO_K$ denotes the ring of integers of $K$. Let $P_K$ denote the set of all prime ideals of $\mcO_K$.
In this section, we study the solutions of the generalized fruit Diophantine equation, namely
\begin{equation}
\label{GFDE}
ax^d-y^2-z^2 +xyz-c=0
\end{equation} 
over $K$, where $d \geq 3$ is an integer and $a,c \in \mcO_K \setminus \{0\}$. 
Let $T_K:= \{ \mfP \in P_K :\ \e(\mfP|2)=1=\f(\mfP|2)\}$, where $\e(\mfP|2)$ and $\f(\mfP|2)$ denote the ramification index and the inertia degree of the prime $\mfP$ lying above $2$, respectively. 
%	Clearly $T_K \subseteq S_K$.

\subsection{Main Results}
\label{section for main result of x^2=By^p+2^rz^p} 
\begin{theorem}
\label{main result1 for GFDE}
Let $K$ be a number field with $T_K \neq \emptyset$. Let $a,b\in \mcO_K \setminus \{0\}$ and $c=2^db-3^r$ with integers $r\geq 2$ and $d \geq 3$ odd. Then the Diophantine equation $ax^d-y^2-z^2 +xyz-c=0$ has no solution $(x_0, y_0, z_0) \in \mcO_K^3$ with $2 | x_0$.
\end{theorem}

\begin{remark}
\label{Remark to thm1}
Note that, if $2$ splits completely in the field $K$, then $T_K \neq \emptyset$. Hence, the conclusion of Theorem~\ref{main result1 for GFDE} holds over all the number fields $K$ in which $2$ splits completely.
\end{remark}

The following corollary gives explicit values of square-free integers $d$ such that the hypothesis of Theorem~\ref{main result1 for GFDE} holds over the quadratic field $K=\Q(\sqrt{d})$.

\begin{corollary}
\label{cor for quadratic field}
Let $t$ be a square-free integer and let $K=\Q(\sqrt{t})$. Then the hypothesis of Theorem~\ref{main result1 for GFDE} holds over $K$ if and only if $t \equiv 1 \pmod 8$.
\end{corollary}
\begin{proof}
Note that, $2$ splits completely in $K=\Q(\sqrt{t})$ if and only if $t \equiv 1 \pmod 8$ (see \cite[Theorem 25]{M18}). By the definition of $T_K$, we conclude that $2$ splits completely in $K=\Q(\sqrt{t})$ if and only if $T_K \neq \emptyset$. Hence, the proof of the corollary follows from Theorem~\ref{main result1 for GFDE} and Remark~\ref{Remark to thm1}.
\end{proof}
Next, we will calculate the density of the set of all square-free integers $t \geq 2$ such that the hypothesis of Theorem~\ref{main result1 for GFDE} holds over $K=\Q(\sqrt{t})$.
Let $$\N^{\text{sf}}:=\{t\in \Z_{\geq 2} : t \text{ is a square-free integer}\}.$$
We shall now define the relative density of any subset $S \subseteq \N^{\text{sf}}$.
\begin{definition}
\label{def rel density}
For any subset $S \subseteq \N^{\text{sf}}$, the \textit{relative density} of $S$ is defined by 
\[ \delta_{\text{rel}}(S):= \lim_{X \to \infty}\frac{\# \{d \in S : d \leq X\}}{\# \{d \in \N^{\text{sf}} : d \leq X\}},\]
if the above limit exists.
\end{definition}
The following theorem computes the density of the set of all square-free integers $t\geq 2$ such that the hypothesis of Theorem~\ref{main result1 for GFDE} holds over $K=\Q(\sqrt{t})$. 
%		More precisely;
\begin{theorem}
\label{thm for density}
Let $U:=\{t \in \N^{\text{sf}} : t \equiv 1 \pmod 8\}$. Then $ \delta_{\text{rel}}(U)= \frac{1}{6}.$ In particular, if $t \in U$, then the hypothesis of Theorem~\ref{main result1 for GFDE} holds over $K=\Q(\sqrt{t})$. 
\end{theorem}
Combining Corollary~\ref{cor for quadratic field} and  Theorem~\ref{thm for density}, we conclude that the set of square-free integers $t \geq 2$ such that the hypothesis of Theorem~\ref{main result1 for GFDE} holds over $K=\Q(\sqrt{t})$ has density exactly $\frac{1}{6}.$

\section{Proofs of Theorem~\ref{main result1 for GFDE} and Theorem~\ref{thm for density}}
\label{steps to prove main results}
\subsection{Proof of Theorem~\ref{main result1 for GFDE}}
To prove this theorem, we need the following lemma. Recall that $P_K$ is the set of all prime ideals of $\mcO_K$ and $T_K= \{ \mfP \in P_K :\ \e(\mfP|2)=1=\f(\mfP|2)\}$.
\begin{lemma}
\label{key lemma}
Let $K$ be a number field. For any prime ideal $\mfP \in T_K$, we have $\mcO_K / {\mfP^n} \simeq \Z/ 2^n\Z$ for all integers $n\geq 1$.
\end{lemma}
\begin{proof}
Let $\mfP \in T_K$. By the definition of $T_K$, we have $\e(\mfP|2)=1=\f(\mfP|2)$. Since $\f(\mfP|2)=1$, we have $\mcO_K / {\mfP} \simeq \Z/ 2\Z$. Since $\mcO_K / {\mfP}\simeq{\mfP^r}/{\mfP^{r+1}}$ for all $r \geq 1$, it follows that $|\mcO_K / {\mfP}|=2=|{\mfP^r}/{\mfP^{r+1}}|$ for all $r \geq 1$, and hence $|\mcO_K / {\mfP^n}|=2^n$. Since $\e(\mfP|2)=1$, we get $v_\mfP(2^{n-1})=n-1$ and therefore $2^{n-1} \notin \mfP^n$. Hence, $2^{n-1}$ is a non-zero element of the quotient ring $\mcO_K / {\mfP^n}$. Since any ring of order $2^n$ in which $2^{n-1} \neq 0$ is isomorphic to $ \Z/ 2^n\Z$, we conclude that $\mcO_K / {\mfP^n} \simeq \Z/ 2^n\Z$.
\end{proof}
\begin{proof}[Proof of Theorem~\ref{main result1 for GFDE}.]
We will prove this theorem by contradiction. Suppose $(\alpha, \beta, \gamma) \in \mcO_K^3$ is a solution of Equation~\eqref{GFDE} with $2 |\alpha$. Then $\alpha= 2 \alpha_1$, for some  $\alpha_1 \in \mcO_K$. So, we have $ a(2\alpha_1)^d-\beta^2-\gamma^2 +2\alpha_1 \beta \gamma-c=0$, which reduces to the following equation: 
$$\beta^2-2\alpha_1 \beta \gamma+ \gamma^2= a (2\alpha_1)^d-c.$$
This gives $(\beta - \alpha_1 \gamma)^2- ({\alpha_1}^2 -1) \gamma^2= 2^d {\alpha_1}^d a-c$. Take $Y=\beta - \alpha_1 \gamma$ and $Z=\gamma$. Then $Y, Z \in \mcO_K$, and we have
\begin{equation}
	\label{red1 eqn}
	Y^2-({\alpha_1}^2 -1)Z^2= 2^d {\alpha_1}^d a-c.
\end{equation}
By assumption $T_K \neq \emptyset$ and choose $\mfP \in T_K$. By Lemma~\ref{key lemma}, we have $\mcO_K / \mfP \simeq \Z/ 2\Z$ and $\mcO_K / {\mfP^2} \simeq \Z/ 4\Z$.
We now consider two cases.

\noindent {\tt Case 1}: $\mfP | \alpha_1$.
This gives $\alpha_1 \equiv 0 \pmod \mfP$, and hence $\alpha_1^2 \equiv 0 \pmod {\mfP^2}$. Since $c=2^db-3^r$ with $d \geq 2$ is an integer, Equation~\eqref{red1 eqn} reduces to the following equation: 
\begin{equation}
	\label{case1 eqn}
	Y^2+Z^2 \equiv 3^r \pmod {\mfP^2}. 
\end{equation}

Since $\mcO_K / \mfP \simeq \Z/ 2\Z$ and $Y \in \mcO_K$, we have either $Y \equiv 0 \pmod \mfP$ or $Y \equiv 1 \pmod \mfP$. Similarly for $Z$, we have either $Z \equiv 0 \pmod \mfP$ or $Z \equiv 1 \pmod \mfP$. If $Y \equiv 0 \pmod \mfP$, then $Y^2 \equiv 0 \pmod {\mfP^2}$. If $Y \equiv 1 \pmod \mfP$, then $(Y-1)^2\equiv 0 \pmod {\mfP^2}$. Since $\mcO_K / {\mfP^2} \simeq \Z/ 4\Z$, we get $Y \equiv 1$ or $3 \pmod {\mfP^2}$. So $(Y-1)^2\equiv Y^2-1 \pmod {\mfP^2}$. Since $(Y-1)^2\equiv 0 \pmod {\mfP^2}$, we get $Y^2 \equiv 1 \pmod {\mfP^2}$. 
In both cases, we have either $Y^2 \equiv 0 \pmod {\mfP^2}$ or $Y^2 \equiv 1 \pmod {\mfP^2}$. 
Similarly, we have either $Z^2 \equiv 0 \pmod {\mfP^2}$ or $Z^2 \equiv 1 \pmod {\mfP^2}$. 

Hence, $Y^2+Z^2 \equiv 0$ or $1$ or $2 \pmod {\mfP^2}$. Since $r$ is odd and $\mcO_K / {\mfP^2} \simeq \Z/ 4\Z$, we have $3^r \equiv 3  \pmod {\mfP^2}$. By Equation~\eqref{case1 eqn}, we get $Y^2+Z^2 \equiv 3 \pmod {\mfP^2}$, which is not possible.

\noindent {\tt Case 2}: $\mfP \nmid \alpha_1$. 
Since $\mcO_K / \mfP \simeq \Z/ 2\Z$, we have $\alpha_1 \equiv 1 \pmod \mfP$. Using the same argument as in the previous case, we get $\alpha_1^2 \equiv 1 \pmod {\mfP^2}$. Since $c=2^db-3^r$ and $d \geq 2$ is an integer, Equation~\eqref{red1 eqn} reduces to the following equation: 
\begin{equation}
	\label{case2 eqn}
	Y^2 \equiv 3^r \pmod {\mfP^2}. 
\end{equation}
Since $r$ is odd and $\mcO_K / {\mfP^2} \simeq \Z/ 4\Z$, we get $3^r \equiv 3   \pmod {\mfP^2}$, and hence $Y^2 \equiv 3 \pmod {\mfP^2}$. This is not possible since $Y^2 \equiv 0$ or $ 1 \pmod {\mfP^2}$.
This completes the proof of the theorem.
\end{proof}

\subsection{Proof of Theorem~\ref{thm for density}}
Recall that $\N^{\text{sf}}=\{t\in \Z_{\geq 2} : t \text{ is a square-free integer}\}.$
Before proving this theorem, we first recall the absolute density of any subset $S \subseteq \N$ (see \cite[Section~7]{FS15}).

\begin{definition}
\label{def for density}
For any subset $S \subseteq \N$ and a positive real number $X$, let $S(X):=\{d \in S : d \leq X\}.$ Then the \textit{absolute density} of $S$ is defined by 
$$ \delta_{\text{abs}}(S):= \lim_{X \to \infty}\frac{\# S(X)}{X},$$
if the above limit exists.
\end{definition}	
The following theorem is useful in the proof of Theorem~\ref{thm for density} (see \cite[Theorem 10]{FS15}).
\begin{theorem}[\cite{FS15}]
\label{FS1}
For $r\in \Z$ and $N \in \N$, let $\N_{r,N}^{\text{sf}}:=\{t\in \N^{\text{sf}} : t \equiv r \pmod N \}$. If $s:= \gcd (r, N)$ is square-free, then 
$\#\N_{r,N}^{\text{sf}} (X) \sim \frac{\varphi(N)}{s\varphi(\frac{N}{s})N \prod_{q|N}(1- \frac{1}{q^2})} \times \frac{6}{\pi^2} X,$
where $\varphi$ denotes Euler's totient function and $q$ varies over all the rational primes dividing~$N$.
\end{theorem}
We are now ready to prove Theorem~\ref{thm for density}.
\begin{proof}[Proof of Theorem~\ref{thm for density}]
Since $\N^{\text{sf}}= \N_{0,1}^{\text{sf}}$, by Theorem~\ref{FS1}, we have $\#\N^{\text{sf}}(X) \sim \frac{6}{\pi^2} X$. By Definition~\ref{def for density}, we have $\delta_{\text{abs}} (\N^{\text{sf}})= \frac{6}{\pi^2}$ (see \cite[page 635]{L09} for more details). Using Definitions~\ref{def rel density} and~\ref{def for density}, we conclude that for any subset $S \subseteq \N^{\text{sf}}$, the absolute density $ \delta_{\text{abs}}(S)$ exists if and only if the relative density $ \delta_{\text{rel}}(S)$ exists. In particular, we have 
\begin{equation}
	\label{rel between density}
	\delta_{\text{abs}}(S) =\delta_{\text{abs}} (\N^{\text{sf}}) \times \delta_{\text{rel}}(S)=\frac{6}{\pi^2} \delta_{\text{rel}}(S).
\end{equation}
By Theorem~\ref{FS1}, we have $ \# \N_{1,8}^{\text{sf}} (X)\sim \frac{1}{\pi^2} X$. Using Definition~\ref{def for density}, we get  $\delta_{\text{abs}}( \N_{1,8}^{\text{sf}})= \frac{1}{\pi^2}$. Finally, by Equation~\eqref{rel between density}, we have $\delta_{\text{rel}}(\N_{1,8}^{\text{sf}})= \frac{\pi ^2}{6} \times  \frac{1}{\pi^2}= \frac{1}{6}$. This completes the proof of the theorem. 
\end{proof}	

\section{Applications}
\label{section for application}
In this section, we will give several applications of the first main result of this article, i.e., Theorem~\ref{main result1 for GFDE}. For the first application, we construct infinitely many elliptic curves $E$ defined over $K$ such that $E$ has no integral point $(x_0,y_0) \in \mcO_K^2$ with $2|x_0$.
%		 More precisely;
\begin{theorem}
\label{appl thm}
Let $K$ be a number field with $T_K \neq \emptyset$ and let $\alpha \in \mcO_K$ be an element not satisfying the polynomial $x^8+5x^6-432x^4-4320x^2-10800$. Let $E_\alpha/K$ be the elliptic curve defined over $K$ given by the Weierstrass equation 
\begin{equation}
	\label{inf elliptic curve}
	E_\alpha: y^2-\alpha xy=x^3-(\alpha^2+5).
\end{equation}
Then $E_\alpha/K$ has no integral point $(x_0,y_0) \in \mcO_K^2$ with $2|x_0$.
\end{theorem} 
\begin{proof}
Since the discriminant $\Delta_{E_\alpha}$ of $E_\alpha$ is $\alpha^8+5\alpha^6-432\alpha^4-4320\alpha^2-10800$ and $\alpha$ is not a root of the polynomial $x^8+5x^6-432x^4-4320x^2-10800$, we get $\Delta_{E_\alpha} \neq 0$. Hence, $E_\alpha$ is an elliptic curve defined over $K$. Now, we will prove this theorem by contradiction.

Suppose $E_\alpha/K$ has an integral point $(x_1,y_1) \in \mcO_K^2$ with $2|x_1$. This gives 
\[y_1^2-\alpha x_1y_1=x_1^3-(\alpha^2+5).\] Hence, $x_1^3-y_1^2-\alpha^2+\alpha x_1y_1-5=0$. Therefore, $(x_1,y_1, \alpha) \in \mcO_K^3$ is an integral solution of Equation~\eqref{GFDE} with $a=b=1$, $r=1$ and $d=3$. This contradicts Theorem~\ref{main result1 for GFDE} since $2|x_1$ and $T_K \neq \emptyset$. Hence, the proof of the theorem follows.
\end{proof}
\begin{remark}
Since the polynomial $x^8+5x^6-432x^4-4320x^2-10800$ has at most $8$ solutions in $K$, the construction of elliptic curves $E_\alpha$ in Equation~\eqref{inf elliptic curve} holds for almost all $\alpha \in \mcO_K$. Note that in \cite{VS22}, Vaishya and Sharma constructed elliptic curves $E_m/ \Q$ for integers $m$ such that $E_m$ has no point $(x_0,y_0) \in \Z^2$, while in Theorem~\ref{appl thm} we construct elliptic curves $E_\alpha$ for almost all $\alpha \in \mcO_K$ such that $E_\alpha/K$ has no integral point $(x_0,y_0) \in \mcO_K^2$ with $2|x_0$.
\end{remark}
To give the next application of Theorem~\ref{main result1 for GFDE}, we first recall the following result (see \cite[Chapter VIII, Theorem 7.1(a)]{S86}).
\begin{theorem}[\cite{S86}]
\label{Sil thm1}
Let $K$ be a number field and let $E/K$ be an elliptic curve given by the  Weierstrass equation
\[E: y^2+a_1xy+a_3y= x^3 + a_2x^2 + a_4x + a_6, \]
where $a_i \in \mcO_K$ for all $i$. If $P=(x,y)\in E(K)$ is a torsion point of order $m \geq 2$ which is not a prime power, then $x,y \in \mcO_K$.
\end{theorem}
As a combination of Theorems~\ref{appl thm} and~\ref{Sil thm1}, we construct infinitely many elliptic curves $E$ defined over $K$ having no torsion point $P\in E(K)$ of order $m \geq 2$ which is not a prime power. 
%		More precisely;
\begin{corollary}
Let $ \alpha, E_\alpha$ be as in Theorem~\ref{appl thm}. Then for each $\alpha$, the elliptic curve $E_\alpha $ has no non-trivial torsion point $P=(x_0,y_0)\in E_\alpha (K)$ which is not a prime power order such that $2|x_0$ and $x_0 \in \mcO_K$.
\end{corollary}
We now recall the Nagell--Lutz theorem (see \cite[Chapter VIII, Corollary 7.2(a)]{S86})).
\begin{theorem}[\cite{S86}]
\label{Sil thm2}
Let $E/ \Q$ be an elliptic curve defined over $\Q$ given by the Weierstrass equation 
\[E: y^2= x^3 + Ax + B, \]
where $A,B \in \Z$. If $P=(x,y)\in E(\Q)$ is a torsion point of order $m \geq 2$, then $x,y \in \Z$.
\end{theorem}
As a combination of Theorems~\ref{appl thm} and~\ref{Sil thm2}, we construct infinitely many elliptic curves $E$ defined over $\Q$ having no torsion point $P\in E(\Q)$ of order $m \geq 2$. 
%		More precisely;
\begin{corollary}
Let $K=\Q$ and let $ \alpha, \ E_\alpha$ be as in Theorem~\ref{appl thm}. Then for each $\alpha$, the elliptic curve $E_\alpha $ has no non-trivial torsion point $P=(x_0,y_0)\in E_\alpha (\Q)$ with $2|\mathrm{Num}(x_0)$, where $\mathrm{Num}(x_0)$ denotes the numerator of the fraction $x_0$ in lowest form.
\end{corollary}
	
		\section*{Acknowledgements}  
%	\vskip20pt\noindent {\bf Acknowledgements.}
	The authors express their sincere gratitude to Professor Kalyan Chakraborty for engaging in discussions regarding the article \cite{PC23} during his visit to the Indian Institute of Technology Hyderabad in September 2023. The authors also extend their sincere gratitude to Professor Cornelius Greither for his insightful mathematical comments towards the improvement of this article. The authors thank Dr. Om Prakash for answering our questions through various emails. The authors are grateful to the editor and the anonymous referee for carefully reading the article and for their insightful suggestions and comments, which have greatly improved the manuscript. The second named author was supported by ANRF CRG Project CRG/2023/005564 while working on this project.

%	\begin{thebibliography}{0}

\end{document}